\documentclass{article}
\usepackage{amsfonts,amsmath,amssymb,url}
\usepackage{hyperref}
\usepackage{breakurl}
\begin{document}
\bibliographystyle{plain}
\begin{title}
{Note on Computing Ratings from Eigenvectors
\thanks{Copyright \copyright\ 2010, R. P. Brent.
\hspace*{\fill}rpb237}
}
\author{Richard P.\ Brent\\
Mathematical Sciences Institute\\
Australian National University\\
Canberra, ACT 0200, Australia}
\date{In memory of\\ Gene Howard Golub\\[2pt] 1932--2007}
\end{title}
\maketitle

\begin{abstract}
We consider the problem of computing ratings using the results of games
played between a set of $n$ players, and show how this problem can be
reduced to computing the positive eigenvectors corresponding to the dominant
eigenvalues of certain $n \times n$ matrices.  There is a close connection
with the stationary probability distributions of certain Markov chains. In
practice, if $n$ is large, then the matrices involved will be sparse, and
the power method may be used to solve the eigenvalue problems efficiently.
We give an algorithm based on the power method, and also derive the
same algorithm by an independent method.
\end{abstract}

\section{Introduction} 
\label{sec:intro}

\thispagestyle{empty}
Suppose that $n$ players, numbered $1,\ldots,n$, play a total of $m$ games
which may end in a win, loss or draw. We assume that a win scores $1$ point,
a draw $0.5$, and a loss $0$. The results can be summarised by an 
$n \times n$ {\em score matrix} $S = (s_{i,j})$, 
where $s_{i,j}$ is the number of points that player $i$ scores against
player $j$ ($i \ne j$). The diagonal entries $s_{i,i}$ are arbitrary; for
reasons discussed later we assume that
$s_{i,i} = \sigma$, where $\sigma \ge 0$ is a constant. For the moment,
the reader may assume that $\sigma = 0$.

The aim is to assign {\em ratings} $r_i$ to the players in such a way that
players who perform better obtain higher ratings. The problem may arise when
$S$ represents the results obtained in a single event, e.g. a Swiss
tournament, and typically we want to use the ratings to break ties between
players on equal scores. It may also arise when $S$ represents all the recent
results involving a large set of players, for example in the regular updates
to a national rating system. In the latter case, we can expect
the matrix $S$ to be large and sparse.
 
The expected score of one player when playing against another
should depend only on the difference of their ratings.
Thus, we assume that, if the
ratings are correct, then the expected score of player $i$
in a game against player $j$ is $f(r_i - r_j)$,
for some function $f: {\cal R} \rightarrow [0,1]$.
Since the total score obtained by both players in a game is $1$,
the function $f$ should satisfy the condition
\begin{equation}
f(z) + f(-z) = 1\;.				\label{eq:fcond1}
\end{equation}
It is reasonable to assume that $f(z)$ is monotonic increasing, and that
\begin{equation}
0 = \lim_{y\to -\infty} f(y) < f(z) < \lim_{y \to +\infty} f(y) = 1\;.
						\label{eq:fcond2}
\end{equation}

It is easy to find functions satisfying these conditions.  For example, if
$\phi(x)$ is a probability density on $(-\infty,+\infty)$,
satisfying the condition $\phi(x) > 0$ and $\phi(x) = \phi(-x)$,
then 
\[f(z) = \int_{-\infty}^z \phi(x)\,dx\]
satisfies~(\ref{eq:fcond1}) and~(\ref{eq:fcond2}).
We could take the normal probability density
\[\phi(x) = \frac{1}{\sqrt{2\pi}}\exp(-x^2/2)\;,\]
but there seems to be no real justification for this choice,
and we shall make a more convenient choice below.  

Let $g(z) = f(z)/(1-f(z))$. For a game without draws, $g(r_i-r_j)$
is the ratio \mbox{$Prob$(player $i$ wins)$/Prob$(player $j$ wins)}.
From~(\ref{eq:fcond1}), $g(z) = f(z)/f(-z)$ and
\begin{equation}
	g(z)g(-z) = 1\;.			\label{eq:gfneqn}
\end{equation}
A simple solution to~(\ref{eq:gfneqn}) is $g(z) = e^{cz}$ for some 
constant~$c$.
Since $f(z) = g(z)/(1+g(z))$, this suggests taking
\begin{equation}
f(z) = \frac{1}{1 + e^{-cz}}\;.			\label{eq:logistic}
\end{equation}
It is easy to check that $f$ 
satisfies~(\ref{eq:fcond1}) and~(\ref{eq:fcond2}) if $c > 0$.

We shall {\em assume} that $f(z)$ has the form~(\ref{eq:logistic}).
This involves an empirical
assumption that could be tested by experiment.  For example, given
three players $\{1,2,3\}$ and ratings $r_1,r_2,r_3$ such that
player $1$ has expected score $f(r_1-r_2)$ against player $2$, and
player $2$ has expected score $f(r_2-r_3)$ against player $3$, 
is it true that player $1$ has expected score $f(r_1-r_3)$ against
player $3$?  The outcome of an experiment to test this hypothesis
might depend on the particular players and what game they are playing.
In practice our assumptions are computationally convenient and we expect
that they will be a reasonable approximation to the truth.

There is some evidence~\cite{Elo1} that the choice~(\ref{eq:logistic}) 
of $f(z)$ (corresponding
to the {\em Logistic} distribution) gives a better approximation to reality,
at least for chess, than the choice based on the normal distribution, as
originally proposed by Elo~\cite{Elo2}. 
Essentially, (\ref{eq:logistic}) is used in the current USCF rating 
system~\cite{Glickman}.
However, in this paper our choice of~(\ref{eq:logistic}) is
made primarily for computational convenience, and because it leads
to an elegant algorithm.

By scaling the ratings, we can assume that $c=1$.  
Thus, in the
following we assume that
\begin{equation}
f(z) = \frac{1}{1 + e^{-z}}\;.			\label{eq:f}
\end{equation}

To avoid the exponential function, it is convenient to define
\[x_i = \exp(r_i)\;.\]
Thus, player $i$ has expected score
\[\frac{1}{1 + x_j/x_i} = \frac{x_i}{x_i + x_j}\]
in a game against player $j$.

The ratings $r_i$ are given by $r_i = \ln x_i$, but we can add
an arbitrary constant $\beta$ to all the $r_i$, since only their differences
are significant (this corresponds to multiplying all the $x_i$ by 
a positive constant $\kappa = \exp(\beta)$).

If $i \ne j$, the total number of games played between player $i$ and player
$j$ is \[g_{i,j} = s_{i,j} + s_{j,i}\;.\]
In the case $i=j$, we use this as a definition of $g_{i,i}$,
that is $g_{i,i} = 2s_{i,i} = 2\sigma$.

Let $W = (w_{i,j})$ be a symmetric matrix of positive {\em weights}
$w_{i,j}$.  Games between players $i$ and $j$ will be weighted in
proportion to $w_{i,j}$.

The {\em actual weighted score} of player $i$ is
\[s_i = \sum_{j=1}^n s_{i,j}w_{i,j}\]
and, given the players' ratings, the {\em expected weighted score}
is 
\[e_i = \sum_{j=1}^n \left(\frac{x_i}{x_i+x_j}\right)w_{i,j}g_{i,j}\;.\]

If the only information available on the players' strengths is the results
encoded in the matrix $S$, then it is reasonable to choose ratings such that
the expected and actual weighted scores of each player are the same,
that is
\[e_i = s_i \;\;{\rm for}\;\; i = 1, \ldots, n\;.\]
Using the definitions of $e_i$, $s_i$ and $g_{i,j}$, this condition is
\begin{equation}
\sum_{j=1}^n \left(\frac{x_i}{x_i+x_j}\right)w_{i,j}(s_{i,j}+s_{j,i}) =
  \sum_{j=1}^n s_{i,j}w_{i,j}
  \;\;{\rm for}\;\; i = 1, \ldots, n \;. 	\label{eq:consist1}
\end{equation}

\section{Choosing Weights to Give a Linear Problem}
\label{sec:weights}

The system of equations~(\ref{eq:consist1}) is nonlinear in the unknowns
$x_i$, and it is not immediately obvious that it has a solution, or how such
a solution should be found.  It is obvious that a solution is not unique,
because if $x = (x_i)$ is one solution then so is $\kappa x$ for any positive
constant $\kappa$.

In order to obtain a linear problem, we choose
\begin{equation}
w_{i,j} = (x_i + x_j)u_{i,j}			\label{eq:weight}
\end{equation}
for some symmetric positive matrix $U$.

With the choice~(\ref{eq:weight}) of weights,
equation~(\ref{eq:consist1}) reduces to
\[\sum_{j=1}^n x_iu_{i,j}(s_{i,j}+s_{j,i}) = 
  \sum_{j=1}^n s_{i,j}(x_i+x_j)u_{i,j}\;.\]
Using the symmetry of $U$, this simplifies to
\begin{equation}
x_i \sum_{j=1}^n a_{j,i} = \sum_{j=1}^n a_{i,j}x_j\;,
							\label{eq:simple}
\end{equation}
where $a_{i,j} = s_{i,j}u_{i,j}$.  The matrix $A = (a_{i,j})$ is a
weighted version of the score matrix $S$, and has
the same sparsity pattern as $S$.

Let \[d_i = \sum_{j=1}^n a_{j,i}\;.\]
$d_i$ can be interpreted as the (weighted) numher of points
lost by player $i$, that is the (weighted) number of points 
scored by player $i$'s opponents in the games they played against $i$.

If a player does not lose any points, then the data is insufficient to 
determine a finite rating for him~-- he is ``infinitely good''.  Thus,
we assume that $d_i > 0$ for $1 \le i \le n$.

Let $D = {\rm diag}(d_i)$ be the diagonal matrix with diagonal elements $d_i$.
By our assumption, $D$ is nonsingular.

\section{The Eigenvalue Problem}
\label{sec:eig}

The condition~(\ref{eq:simple}) can be written in matrix-vector form as
\begin{equation}
Dx = Ax					\label{eq:eig1}
\end{equation}
or
\begin{equation}
D^{-1}Ax = x\;.					\label{eq:eig2}
\end{equation}
Thus, the solution vector $x$ is the eigenvector corresponding to the
eigenvalue $1$ of the matrix $D^{-1}A$.

We observe that the matrix $A-D$ has linearly dependent rows; in fact,
it is easy to see from the definition of $D$ that the rows of $A-D$ sum
to zero.  Thus, $A-D$ is singular, so $D^{-1}A - I$ is singular,
and $D^{-1}A$ does in fact have an eigenvalue $1$.  Similarly
for $AD^{-1}$.

Equation~(\ref{eq:eig1}) can be interpreted in terms of a Markov chain.
Let $y = Dx$ and $M^T = AD^{-1}$.  Then $M$ is the transition matrix
of a Markov chain ($m_{i,j} \ge 0$ and $\sum_j m_{i,j} = 1$).
The vector $y/||y||_1$ gives the stationary probability
distribution, because $y^T M = y^T$, or equivalently $AD^{-1}y = y$.
It follows from standard theory of Markov matrices that 
$\rho(D^{-1}A) = \rho(AD^{-1}) \le 1$.

In certain degenerate cases we can not expect finite ratings to
be defined by the data.  We already assumed that $d_i > 0$.
This is necessary, but not sufficient.  If the players
can be split into two disjoint nonempty sets such that
players in the first set always beat the players in the second set, 
then the players in
the first set are ``infinitely better'' than the players in the second set.  
Similarly, if players in the first set never play players in 
the second set, we can not
expect to compare their playing strengths. In practice, in either of these
situations, we could split the problem and rate players in each set separately.  

In the typical nondegenerate case, $D^{-1}A$ has a simple eigenvalue 
$\lambda_1 = 1$,
and the other eigenvalues $\lambda_i$ are inside the unit circle,
that is $|\lambda_i| < 1$ except for $\lambda_1 = 1$.
Then the {\em power method} converges and we can find $x$ by the simple 
iteration
\[x^{(k+1)} = D^{-1}Ax^{(k)}\;,\]
with a suitable starting vector, e.g. $x^{(0)} = (1,1,\ldots,1)^T$.

So far we did not mention the role of the constant $\sigma$ (recall that
$s_{i,i} = \sigma$).  The solution $x$ of~(\ref{eq:eig2}) is independent of
$\sigma$, but the speed of convergence of the power method depends on
$\sigma$. We have found in our experiments that $\sigma \in [0.2, 0.5]$
is a good choice to maximise the speed of convergence.  Any $\sigma > 0$
will ensure that $D$ is nonsingular.

\section{Modifying the Weights}
\label{sec:mod}

We have seen that solving an eigenvalue problem allows us to compute
ratings if the score matrix is weighted by the weight
function~(\ref{eq:weight}). It would be more natural to solve the problem
with unit weights, that is $w_{i,j} = 1$.  Unit weights have the advantage
that, when applied to a round-robin (``all play all'') tournament, players
with the same score obtain the same ratings, as is easy to see
from~(\ref{eq:consist1}).

The condition $w_{i,j} = 1$ is equivalent to
\[u_{i,j} = \frac{1}{x_i + x_j}\;.\]
Thus, we can regard the solution of the eigenvalue problem~(\ref{eq:eig2})
as an inner iteration, and introduce an outer iteration where we change
the weights. If $x^{(k)}$ is the solution to the $k$-th eigenvalue
problem (with $x^{(0)} = (1,1,\ldots,1)^T$ an initial vector),
then the $(k+1)$-st eigenvalue problem will use
\[u_{i,j}^{(k+1)} = \frac{1}{x_i^{(k)} + x_j^{(k)}}\;.\]
If the outer iteration converges, then it solves the original
problem with weights
\[w_{i,j} = (x_i + x_j)u_{i,j} = 1\;.\]

In practice, we have found that convergence is quite rapid in
nondegenerate cases.  However, it is wasteful to solve the inner
eigenvalue problems accurately.  It is much more efficient to perform
just one iteration of the power method in the inner loop.
The resulting Algorithm~1 (without improvements to take advantage of
sparsity) is given below.

\begin{table}[ht]
\label{tab:alg1}
\begin{verbatim}
          for i := 1..n do 
            x[i] := 1.0;
            end for;
          for k := 1, 2, ... until convergence do
            for i := 1..n do
              d[i] := 0.0;
              end for;
            for i := 1..n do
              sum := 0.0; 
              for j := 1..n do
                temp := s[i,j]/(x[i] + x[j]);
                sum := sum + temp*x[j];
                d[j] := d[j] + temp;
                end for;
              y[i] := sum;
              end for;
            for i := 1..n do
              x[i] := y[i]/d[i];
              end for;
            end for;
\end{verbatim}
\begin{center}
{\bf Algorithm 1: Unit Weights}
\end{center}
\end{table}

Since the aim is to compute ratings $r_i = \ln x_i$,
the convergence test should ensure a small relative error in each
component of $x$.  Thus, an appropriate stopping criterion is
\[\max_{1 \le i \le n}\left|\frac{x_i^{(k)}-x_i^{(k-1)}}{x_i^{(k)}}\right| 
  < \varepsilon\;,\]
where $\varepsilon$ is a tolerance depending on the accuracy required. 

Failure to converge in a reasonable number of iterations may
indicate that the problem is degenerate and that some ratings
are diverging to $\pm \infty$.  In this case one or more players should
in principal be excluded from consideration.  A more convenient solution
in practice is to add a ``dummy'' player who draws with all the other
players, and whose games are given a positive weight $\gamma$,
for example $\gamma = 1$. As $\gamma \to 0+$ the ratings tend to
the correct values (possibly $\pm\infty$), but for any positive $\gamma$
we obtain a nondegenerate problem and finite ratings\footnote{This 
solution is similar to the one adopted in the {\em PageRank} 
algorithm used by the Google search engine~\cite{BP,PB,PR2}, 
where a fictional page
essentially has links to every other page. Our parameter $\gamma$
corresponds to the Page Rank algorithm's $1-d$.}.

If $n$ is large then $S$ (and hence $A$) will be sparse, since there are
at most two off-diagonal entries for each game played. 
Thus, the number of nonzero elements is at most $2m + n$.  The inner loop of
the algorithm above essentially involves the multiplication of $A$ on the
right by $x$, and on the left by $[1,1,\ldots,1]$. Thus, standard sparse
matrix techniques can be used to reduce the complexity of the inner
iteration from $O(n^2)$ to $O(m)$.

In practice the final ratings would be modified by a linear transformation to 
make them positive and not too small, before rounding to the nearest integer. 
(FIDE Elo ratings are usually in the range $[0,3000]$,
and BCF ratings are usually in the range $[0,300]$, see~\cite{conversion}.)

\section{An Independent Derivation of Algorithm 1}
\label{sec:indep}

{From} (\ref{eq:consist1}) with $w_{i,j} = 1$ we have, for $i = 1, \ldots, n$,

\begin{equation}
\sum_{j=1}^n \left(\frac{x_i}{x_i+x_j}\right)(s_{i,j}+s_{j,i}) =
  \sum_{j=1}^n s_{i,j}\;.				\label{eq:consist2}
\end{equation}
This simplifies to
\[
\sum_{j=1}^n \frac{x_i s_{j,i}}{x_i+x_j} = 
\sum_{j=1}^n \frac{x_j s_{i,j}}{x_i+x_j}
\]
and, taking $x_i$ outside the sum on the left, we see that
\begin{equation}
x_i = \left(\sum_{j=1}^n \frac{s_{i,j}x_j}{x_i+x_j}\right)\left/
      \left(\sum_{j=1}^n \frac{s_{j,i}}{x_i+x_j}\right)\right.\;.
							\label{eq:consist3}
\end{equation}
A natural iteration to solve~(\ref{eq:consist3}) is
\begin{equation}                                                 
x_i^{(k+1)} = \left(\sum_{j=1}^n \frac{s_{i,j}x_j^{(k)}}
				       {x_i^{(k)}+x_j^{(k)}}\right)\left/ 
      \left(\sum_{j=1}^n \frac{s_{j,i}}{x_i^{(k)}+x_j^{(k)}}\right)\right.
                                                        \label{eq:iter1}
\end{equation}            
for $k = 1, 2, \ldots$.                              
However, it is easy to see that this is exactly the iteration
implemented in Algorithm~1!

The significance of the diagonal terms $s_{i,i} = \sigma \ge 0$ is apparent
if we consider the diagonal terms ($j=i$) in the numerator and denominator
of~(\ref{eq:iter1}). The diagonal term in the numerator is $\sigma/2$
and the diagonal term in the denominator is $\sigma/(2x_i^{(k)})$.
As $\sigma \to \infty$ the right-hand side of~(\ref{eq:iter1})
$\to x_i^{(k)}$.
Thus, $\sigma$ acts as a damping factor: 
larger values of $\sigma$ tend to reduce the change 
$x_i^{(k)} \rightarrow x_i^{(k+1)}$ at each iteration.

Other iterations can be obtained in a similar manner.
For example, taking $x_i$ outside the sum on the left side
of~(\ref{eq:consist2}), we obtain the iteration
\begin{equation}
x_i^{(k+1)} = \left(\sum_{j=1}^n s_{i,j}\right)\left/ 
      \left(\sum_{j=1}^n \frac{s_{i,j}+s_{j,i}}
	{x_i^{(k)}+x_j^{(k)}}\right)\right.\;.
\label{eq:iter2}
\end{equation}
However, our numerical experiments suggest that~(\ref{eq:iter2}) gives
slower convergence than~(\ref{eq:iter1}).  This conclusion is confirmed
by a first-order analysis in the special case that all the $x_i$ are
approximately equal.
 
\section{Incorporating Old Ratings}

Often some or all of the players will already have ratings,
say player $i$ has rating ${\widehat x_i}$ based on ${\widehat w_i}$ 
games\footnote{The weight ${\widehat w_i}$ associated with an old rating
might be reduced by a constant factor, say $0.5$, to give less weight
to old games than to recent ones.}.
It is easy to take such ``old'' ratings into account by a slight
modification of the argument leading to equation~(\ref{eq:consist1}).
We add ${\widehat w_i}/2$ to the actual weighted score $s_i$, 
as if player $i$ drew ${\widehat w_i}$ games against a player with
rating ${\widehat x_i}$, 
and also add ${\widehat w_i}x_i/(x_i + {\widehat x_i})$
to the expected weighted score $e_i$.  Thus equation~(\ref{eq:consist1})
becomes
\begin{equation}
\sum_{j=1}^n \left(\frac{x_i}{x_i+x_j}\right)w_{i,j}(s_{i,j}+s_{j,i}) 
  \;+\; \frac{{\widehat w_i}x_i}{x_i+{\widehat x_i}} \;=\;
  \sum_{j=1}^n s_{i,j}w_{i,j}
  \;+\; \frac{\widehat w_i}{2} 			\label{eq:consist1a}
\end{equation}
for $i = 1, \ldots, n$,
and equation~(\ref{eq:consist2}) becomes 
\begin{equation}
\sum_{j=1}^n \left(\frac{x_i}{x_i+x_j}\right)(s_{i,j}+s_{j,i}) \;=\;
  \sum_{j=1}^n s_{i,j}
  \;-\; \left(\frac{x_i - \widehat x_i}{x_i + \widehat x_i}\right)
  \frac{\widehat w_i}{2}\;.				\label{eq:consist2a}
\end{equation}
Now, it is easy to see that the iteration~(\ref{eq:iter1})
becomes
\begin{equation}                                                 
x_i^{(k+1)} \;=\; \frac{\displaystyle 
	\left(\sum_{j=1}^n \frac{s_{i,j}x_j^{(k)}}
				       {x_i^{(k)}+x_j^{(k)}}\right)
 \;+\; \frac{{\widehat w_i}{\widehat x_i}}{2(x_i^{(k)} + {\widehat x_i})}}
     {\displaystyle 
	\left(\sum_{j=1}^n \frac{s_{j,i}}{x_i^{(k)}+x_j^{(k)}}\right)
 \;+\; \frac{\widehat w_i}{2(x_i^{(k)} + {\widehat x_i})}}
                                                        \label{eq:iter1a}
\end{equation}            
for $k = 1, 2, \ldots$.                              

\pagebreak[3]
\section{Conclusion}
\label{sec:conc}

We have shown how the computation of ratings for players of chess
(and other games) can be reduced to an eigenvalue problem, or a sequence
of eigenvalue problems, and that these eigenvalue problems are closely
related to the problem of computing the stationary probability distributions
of certain Markov chains.  The eigenvalue problems can be solved efficiently
by the power method. We derived an algorithm (Algorithm~1) by two 
different methods; one derivation (\S\S\ref{sec:eig}--\ref{sec:mod})
used the power method, the other (\S\ref{sec:indep}) was from first
principles.

We have tested Algorithm~1 on small examples and on some larger examples
with simulated scores.  It would be interesting to test the algorithm on
real data and to see how it compares with the methods currently in
use for chess~\cite{conversion,Elo1,Elo2,Glickman} and other games,
e.g.~football~\cite{Colley}.  Because they are practical
systems that have evolved over time, most of these methods involve various
{\em ad hoc} features and piecewise linear approximations, so they are less
elegant (though perhaps more practical) than the relatively simple algorithm
discussed here.

\section*{Acknowledgement}

Thanks to my friend and mentor Gene Golub who persuaded me to write up this
work, and commented on an early draft.
Unfortunately, Gene passed away in 2007.


\begin{thebibliography}{11}

\bibitem{Elo1}
Anonymous, ELO rating system,
\url{http://en.wikipedia.org/wiki/Elo\_rating\_system}

\bibitem{BP}
Sergey Brin and Lawrence Page, %
The Anatomy of a Large-Scale
Hyper\-textual Web Search Engine,
{\em Proc.\ Seventh International Web Conference} (WWW~98), 1998.
Also \url{http://www-db.stanford.edu/~backrub/google.html}

\bibitem{Callaghan}
Thomas Callaghan, Peter J.~Mucha and Mason A.~Porter,
The Bowl Championship Series: A mathematical review,
{\em AMS Notices} {\bf{51}}, 8 (2004), 887--893.	%

\bibitem{Colley}
Wesley N. Colley, 
Colley's bias free college football ranking method: the
Colley matrix explained, 2002.
\url{http://www.colleyrankings.com/#method}

\bibitem{Connor}
G. R. Connor and C. P. Grant,
An extension of Zermelo's model for ranking by paired comparisons,
{\em European J.\ Appl.\ Math.\ }{\bf{11}} (2000), 225--247.


\bibitem{Elo2}
Arpad Elo, The Ratings of Chess Players, Past and Present, Arco Pub., 1978.

\bibitem{Glickman}
Mark Glickman, Mark Glickman's ratings page,
\url{http://www.glicko.net/ratings.html}

\bibitem{Keener}
J. P. Keener, The Perron-Frobenius theorem and the ranking of football
teams, {\em SIAM Review} {\bf{35}} (1993), 80--93.

\bibitem{conversion}
Chris Majer,
How grading works and other information,
British Chess Federation, 2002.

\bibitem{PB}
Larry Page, Sergey Brin, Rajeev Motwani, and Terry Winograd,
The PageRank citation ranking: bringing order to the web, 
Technical report, Stanford University Database Group, 1999.
\url{http://dbpubs.stanford.edu:8090/pub/1999-66}

\bibitem{PR2}
Ian Rogers, PageRank explained: 
the Google PageRank algorithm and how it works, 2002.
\url{http://www.iprcom.com/papers/pagerank/}

\end{thebibliography}
\end{document}